\documentclass[11pt]{article}
\usepackage[T1]{fontenc}
\usepackage{lmodern}
\usepackage[margin=1in]{geometry}
\usepackage{amsmath,amssymb,amsthm,mathtools,mathrsfs}
\usepackage{microtype}
\usepackage{tikz}
\usepackage{iftex}
\usepackage{etoolbox}
\usepackage[hidelinks]{hyperref}
\usepackage[color]{changebar}
\usepackage{multirow}
\usepackage{comment}
\usepackage{orcidlink}

\usepackage{hyperref}
\usepackage{doi}

\allowdisplaybreaks[2]
\newtheorem{theorem}{Theorem}
\newtheorem{lemma}{Lemma}[section]

\newtheorem{corollary}{Corollary}[section]

\newcommand{\FF}{\mathbb F}
\newcommand{\KK}{\mathbb K}
\newcommand{\VV}{\mathcal V}
\newcommand{\UU}{\mathcal U}
\newcommand{\BB}{\mathcal B}
\newcommand{\Aux}{\mathcal A}
\newcommand{\SSS}{\mathcal S}
\newcommand{\PP}{\mathcal P}
\newcommand{\DD}{\mathscr D}

\newcommand{\iacr}{\operatorname{iacr}}
\newcommand{\acr}{\operatorname{acr}}
\newcommand{\iocr}{\operatorname{iocr}}
\newcommand{\ocr}{\operatorname{ocr}}
\newcommand{\pcr}{\operatorname{pcr}}
\newcommand{\cro}{\operatorname{cr}}

\makeatletter
\newif\if@inbar
\@inbarfalse
\newcommand{\@currentbarcolor}{none}

\newcommand{\startmarkcolor}[1]{%
  \if@inbar
    \def\@tempa{#1}%
    \ifx\@tempa\@currentbarcolor\relax
    \else
      \cbend
      \cbcolor{#1}%
      \def\@currentbarcolor{#1}%
      \cbstart
    \fi
  \else
    \cbcolor{#1}%
    \def\@currentbarcolor{#1}%
    \cbstart
    \@inbartrue
  \fi
}

\newcommand{\stopmark}{%
  \if@inbar
    \cbend
    \@inbarfalse
    \def\@currentbarcolor{none}%
  \fi
}

\AtEndDocument{\stopmark}
\makeatother

\title{Crossing Numbers of Complete Graphs:\\ the independent odd crossing number}

\author{%
  Marcus Schaefer%
  \thanks{DePaul University, \texttt{mschaefer@cdm.depaul.edu}}
  \phantom{}
  \orcidlink{0000-0001-7005-8599}\and
    Daniel \v{S}tefankovi\v{c}%
  \thanks{University of Rochester, \texttt{stefanko@cs.rochester.edu}}%
\phantom{}  
  \orcidlink{0000-0002-4849-7955}
}

\date{November 8, 2026}
\begin{document}
\maketitle
\vspace{-1.4em}
\begin{abstract}
We show that the Harary--Hill bound on the crossing number of complete graphs as well as the Zarankiewwicz bound on the crossing number of complete bipartite graphs, both recently established by OpenAI, also hold for the independent odd crossing number. This implies that the same bounds hold for several other crossing number variants such as the (independent) algebraic crossing number and the (independent) pair crossing number.  
\end{abstract}

\section{Introduction}

The {\em crossing number} of a  graph is the smallest number of crossings achievable in a drawing of the graph. Crossing number research for many decades has been inspired by two conjectures: Zarankiewicz~\cite{Z53} in the 1950s thought he had proven that
\[\cro(K_{m,n}) = Z(m,n),\]
where 
\[
 Z(m,n)=
 \Big\lfloor\frac{m}{2}\Big\rfloor
 \Big\lfloor\frac{m-1}{2}\Big\rfloor
 \Big\lfloor\frac{n}{2}\Big\rfloor
 \Big\lfloor\frac{n-1}{2}\Big\rfloor.
\]
There was a gap in his proof so the conjecture that his bound is correct is now known as the {\em Zarankiewicz conjecture} (for historical background, see Sz\'{e}kely's~\cite{S16} paper). In the 1960s, Harary and Hill~\cite{HH63} conjectured that the crossing number of the complete graph is given by
\[\cro(K_n) = Z(n),\]
where 
\[
 Z(n)=\frac14
 \Big\lfloor\frac n2\Big\rfloor
 \Big\lfloor\frac{n-1}{2}\Big\rfloor
 \Big\lfloor\frac{n-2}{2}\Big\rfloor
 \Big\lfloor\frac{n-3}{2}\Big\rfloor.
\]

Both conjectures were recently proved by OpenAI~\cite{openai,openaiBipartite2026} and the Harary--Hill bound has already been rewritten for human readers by  crossing number experts: Martin Balko~\cite{balko} and Gelasio Salazar~\cite{salazar}. Balko showed that OpenAI's result holds for the crossing number variant $\ocr_+$ (which we will define below), while Salazar, taking a combinatorial point of view, made progress towards Archdeacon's conjecture (for background on the conjecture see~\cite{AMcQRS17}).

We give a third exposition of the OpenAI proof. Our goal is to strengthen the result by showing that the lower bound also applies to other ways of counting crossings in a graph. We focus on $\iocr(G)$,  the {\em independent odd crossing number} of $G$, which is the smallest number of pairs of independent edges that cross an odd number of times in any drawing of $G$. By definition, $\iocr(G) \leq \cro(G)$.

\begin{theorem}\label{thm:main}
For all $n\ge1$,
\[
\iocr(K_n)=\cro(K_n)=Z(n).
\]
\end{theorem}

The theorem also settles the value for $K_n$ of many other crossing number variants, including the algebraic crossing number, the pair crossing number, and the independent crossing number~\cite[Open Question 2]{S13}. We will formally define these crossing numbers in the next section.  

The result on $\iocr(K_n)$ is non-trivial in the sense that we know that there are graphs for which $\iocr(G) < \cro(G)$~\cite{fulek}. For $\ocr_+$ this question remains open: for all we know $\ocr_+(G) = \cro(G)$, though we do not expect this to be the case. 

The same proof technique 
also yields the independent odd crossing number of $K_{m,n}$. The original OpenAI proof required some more
conceptual machinery in this case.
\begin{theorem}\label{thm:bipartite}
For all $m,n\geq1$,
\[
 \iocr(K_{m,n})=\cro(K_{m,n})=Z(m,n).
\]
\end{theorem}

\quad\\

\paragraph{Approach}
The approach follows the framework used in OpenAI's breakthrough. Each vertex in the drawing is assigned an independent algebraically independent variable 
in a field extension of $\FF_2$. The drawing of $K_n$ is  then connected, using the concept of intersection pairings, to the evaluation of bounded-degree four-variable 
polynomials on quadruples of elements corresponding to pairs of independent edges crossing an odd number of times. Using properties of the drawing and the intersection pairing, 
it is shown that a non-zero polynomial in the space cannot ever evaluate to zero on all the points---establishing a lower bound on the number of pairs of independent edges crossing oddly in terms of the dimension of the space of the polynomials. The bound is within a factor of two of the desired bound and is then boosted by another refined dimension 
argument. An important ingredient, that is also present in the previous proofs, is a ``localization argument'' that gives regularity of a map arising from a local intersection 
pattern around a single vertex.  For the independent odd crossing number the drawing allows dependent edges to intersect, resulting in more general intersection pairings and requiring a much 
stronger localization argument.

\paragraph{AI statement} This work has been produced with significant help of ChatGPT 5.6 Pro (understanding the previous results, exploring proof strategies, proving statements, and consulting alternative exposition paths). The authors  take full responsibility for the presentation and validity of the results.

\section{Crossing Numbers}

We restrict ourselves to drawings of a graph in which each pair of edges crosses at most a finite number of times, no three edges cross in a point, and no edge passes through a vertex. For two edges $e$, $f$ in a drawing $D$ we define $\cro_D(e,f)$ as the number of times $e$ and $f$ cross in $D$. The {\em crossing number of a drawing} $D$ is defined as
\[\cro(D) = \sum_{\{e,f\}} \cro_D(e,f).\]
With that the {\em crossing number} of $G$ is the smallest $\cro(D)$ over all drawings $D$ of $G$. 
%
If we count crossings modulo $2$ we obtain the {\em odd crossing number} of $D$:
\[\ocr(D) = \#\{\{e,f\}: \cro_D(e,f) = 1\bmod{2}\}.\]
The {\em odd crossing number} of $G$, $\ocr(G)$, minimizes $\ocr(D)$ over all drawings $D$ of $G$. If we only want to count the crossing pairs in a drawing, we set
\[\pcr(D) = \#\{\{e,f\}: \cro_D(e,f) \geq 1\},\]
and define the corresponding {\em pair crossing number}, $\pcr(G)$ as the smallest $\pcr(D)$ where $D$ ranges over all drawings of $G$. 

We can also count crossings algebraically, an approach we will take in the proof (before reducing to the modulo $2$ version). For two directed edges $e$, $f$ in a drawing $D$
we label each crossing $+1$ or $-1$ depending on the direction in which $e$ crosses $f$. The {\em signed algebraic intersection number} $I_D(e,f)$ is the sum of the signed crossings of $e$ and $f$; flipping the direction of an edge only changes the sign
of $I_D(e,f)$, so we can define the {\em algebraic crossing number} of $D$ as
\[\acr(D) = \sum_{\{e,f\}} |I_D(e,f)|.\]
The {\em algebraic crossing number} of $G$, $\acr(G)$, is the minimum of $\acr(D)$ over all drawings of $G$.
By definition, $\ocr(G) \leq \acr(G) \leq \cro(G)$ for all graphs $G$. 

Pach and T\'{o}th~\cite{PT00} initiated a systematic study of crossing number variants by introducing  two rules that can be used with any crossing number: ``Rule $+$'' prohibits crossings between adjacent edges, while ``Rule~$-$'' allows crossings between adjacent edges, but does not count them. Applying these rules to $\ocr$, $\acr$ and $\cro$ yields eight variants, since $\cro_+ = \cro$ (there are no crossings between adjacent edges in a $\cro$-minimal drawing).

The tables below summarize these variants, they are based on figures from~\cite{BMP05,fulek}. 

\begin{center}
\setlength{\tabcolsep}{4pt}
\renewcommand{\arraystretch}{1.2}
\begin{tabular}{l|l|l}
 $\ocr_+$ & $\acr_+$ & \multirow{2}{*}{$\cro$} \\ \cline{1-2}
 $\ocr$ & $\acr$ & \\ \hline
 $\ocr_-=\iocr$ & $\acr_-=\iacr$ & $\cro_-$
\end{tabular}
\qquad\begin{tabular}{l|l|l}
 $\ocr_+$ & $\pcr_+$ & \multirow{2}{*}{$\cro$} \\ \cline{1-2}
 $\ocr$ & $\pcr$ & \\ \hline
 $\ocr_-=\iocr$ & $\pcr_-$ & $\cro_-$
\end{tabular}
\end{center}


The values of the crossing number variants in each table are non-decreasing as one moves rightwards and/or upwards. Since $\iocr$ and $\cro$ are at extreme ends of both tables, Theorem~\ref{thm:main} implies that the Harary--Hill bound holds for all of these crossing number variants.\footnote{The exact relationship between $\acr$ and $\pcr$ remains open, mainly because there are no known examples for $\pcr<\cro$.} 

Balko~\cite{balko} proves the Harary--Hill bound for $\ocr_+$, and thereby $\acr_+$ (and $\pcr_+$), so his proof can handle the odd crossing count, but only if adjacent crossings are not allowed. Our main contribution is to show that not only can we allow adjacent crossings when counting crossings modulo $2$, but that the crossing number does not change even if we do not count those  crossings. 

The main difference from the previous proofs and expositions is allowing dependent edges to cross. This changes 
how intersection pairing is defined (the second line in~\eqref{eq:Rdirect} accounts for such intersections) and
this, in turn, requires a stronger claim for the localization of the polynomial method to a vertex (Theorem~\ref{thm:detector}). 

For the proof we will need the following terminology:
The {\em rotation} at a vertex $v$ is a cyclic permutation of the edges incident to $v$; a {\em rotation system} of a graph assigns a rotation to each vertex of the graph. A drawing of a graph determines a rotation system by letting the rotation at a vertex $v$ be the counterclockwise permutation of the edges incident to $v$ as they attach to $v$.

\section{Edge Chains and Intersection Pairings}\label{sec:geometry}

Our proof, just like the other proofs~\cite{openai,balko,salazar}, will be based on edge-chains, an idea introduced into crossing number research by Tutte in his famous ``Toward a theory of crossing numbers''~\cite{T70}.

Write $[uv]=-[vu]$ for the directed edge from $u$ to $v$, and $\partial[uv]=[v]-[u]$. Let $[uu]=0$ (note that this convention is compatible with the $\partial$ 
notation: $\partial [uu] = [u]-[u]=0$). An {\em (edge) chain} is an arbitrary linear combination of these 
symbols; a cycle has zero boundary (paths and cycles are special cases of an edge chain). 

Our goal in this section is to explain how $I_D$ can be extended to chains in $G$. To this end we want to lift $I_D$ to pairs of curves $\gamma,\eta$. This is easy if $\gamma$ and $\eta$ do not have any {\em interfering endpoints}, that is if $\partial\gamma \cap \eta = \emptyset$ 
and $\partial\eta\cap\gamma = \emptyset$ (for example, a pair of closed curves has no interfering endpoints):
 For two oriented curves 
with no interfering endpoints the signed algebraic intersection number $I_D(\gamma,\eta)$ is well-defined (the value is robust under 
sufficiently small perturbations of the curves, even if they share internal vertices or segments, or cross).  

To carefully account for intersections of edges in the drawing we will need to extend the signed algebraic intersection number to curves that have interfering endpoints.
We need an {\em intersection pairing}---a bilinear function $R'_D(z,z')$ where 
$z,z'$ are edge chains such that $R'_D$ agrees with $I_D$ for chains $\gamma,\eta$ with no interfering endpoints. The choice of function $R'_D$ is not unique; we
will use a particularly simple form. For every vertex $v$ we take the cyclic permutation of the edges around $v$ and turn it into a linear order $\prec_v$ by 
arbitrarily selecting one of the edges as the first edge in the rotation. 
The function 
$R'_D$ can then be defined as the bilinear extension of the following equations:
\begin{equation}\label{eq:Rdirect}
\boxed{\begin{aligned}
 R'_D([uv],[wq])&=I_D([uv],[wq])
    &&\text{if the edges are independent},\\
 R'_D([vp],[vq])&=I_D([vp],[vq])+\mathbf1\{p\prec_v q\}
    &&\text{if }p\ne q,\\
 R'_D([uv],[uv])&=1.&&
\end{aligned}}
\end{equation}

The following lemma justifies this definition of $R'_D$.
\begin{lemma}\label{lem:extension}
The bilinear extension of~\eqref{eq:Rdirect} agrees with the signed algebraic intersection number $I_D$ for chains $\gamma,\eta$ with non-interfering endpoints.
\end{lemma}

\begin{figure}[t]
\centering
\begin{tikzpicture}[scale=.9,>=stealth]
  \definecolor{redpath}{RGB}{180,54,58}
  \definecolor{bluepath}{RGB}{35,102,168}
  \begin{scope}[shift={(0,0)}]
    \draw[gray!65] (0,0) circle (1.55);
    \coordinate (rh) at (110:1.55);
    \coordinate (bh) at (128:1.55);
    \coordinate (re) at (205:1.55);
    \coordinate (be) at (220:1.55);
    \coordinate (rf) at (300:1.55);
    \coordinate (bf) at (315:1.55);
    \draw[redpath,very thick,->] (rh)--(re);
    \draw[bluepath,very thick,->] (bh)--(bf);
    \foreach \a in {rh,re,rf} \fill[redpath] (\a) circle (2.0pt);
    \foreach \a in {bh,be,bf} \fill[bluepath] (\a) circle (2.0pt);
    \node[font=\scriptsize,above left] at (rh) {$r_v$};
    \node[font=\scriptsize,above right] at (bh) {$b_v$};
    \node[font=\scriptsize,below left] at (re) {$r_e$};
    \node[font=\scriptsize,below left] at (be) {$b_e$};
    \node[font=\scriptsize,right] at (rf) {$r_f$};
    \node[font=\scriptsize,right] at (bf) {$b_f$};
    \node[font=\small] at (0,-2.05) {(a) $e\prec_v f$: local crossing $+1$};
  \end{scope}
  \begin{scope}[shift={(5.3,0)}]
    \draw[gray!50,rounded corners=2pt] (-1.75,-.52) rectangle (1.75,.52);
    \draw[redpath,very thick,->] (-1.65,-.28)--(1.65,.28);
    \draw[bluepath,very thick,->] (-1.65,.28)--(1.65,-.28);
    \node[font=\scriptsize,left,redpath] at (-1.65,-.28) {red};
    \node[font=\scriptsize,left,bluepath] at (-1.65,.28) {blue};
    \node[font=\small,align=center] at (0,-2.05) {(b) Same-edge corridor:\\crossing $-1$};
  \end{scope}
\end{tikzpicture}
\caption{The ordered red--blue resolution.  Only the red copy of $e$
 and blue copy of $f$ are drawn in (a); the other ports illustrate
 their order.  The two copies of a single edge switch sides in (b).}
\label{fig:redblue}
\end{figure}
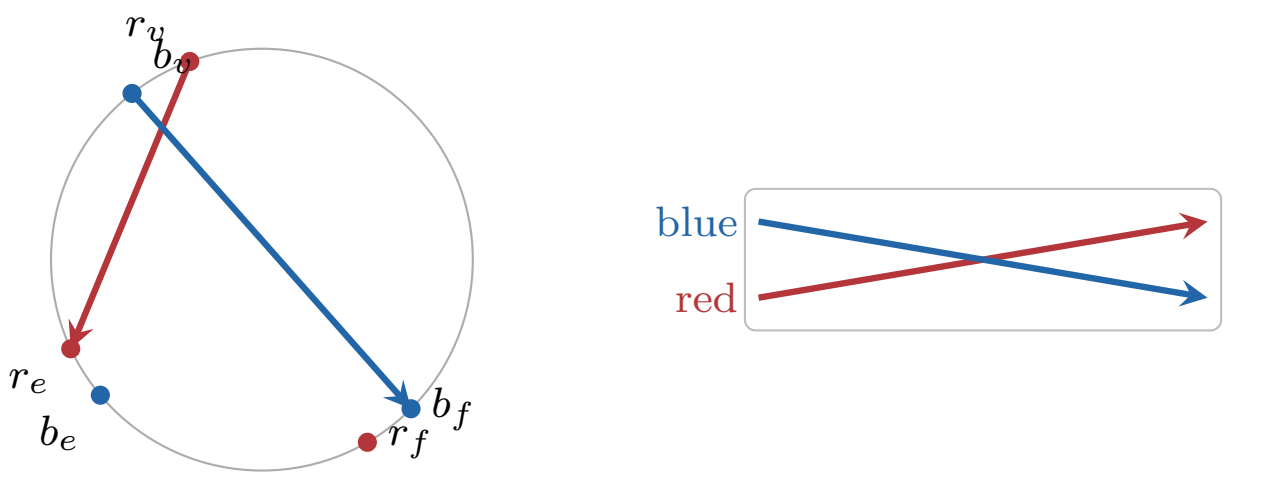

\begin{proof}
For the proof we construct the intersection pairing by a systematic perturbation of the drawing $D$ such that for any two curves whose signed algebraic intersection number we are computing these curves do not have interfering endpoints; the bilinearity then follows
from the bilinearity of the signed algebraic intersection number. 

Take a small disk around each vertex $v$, meeting its incident edges in
spokes.  Choose a gap between two spokes at $v$ declaring them the first and last edges: $e_1\prec_v\cdots\prec_v e_d$.  On the disk boundary add a red vertex
$r_v$ and a blue vertex $b_v$ within the gap, and a red and a blue vertex, which we call {\em ports}, beside
each spoke, the red port slightly clockwise of the blue one.  Thus the
boundary points occur counterclockwise in the order
\begin{equation}\label{eq:ports}
 r_v,\ b_v,\ r_{v,e_1},\ b_{v,e_1},\ldots,
 r_{v,e_d},\ b_{v,e_d}.
\end{equation}
Join each colored vertex to its ports by straight chords inside the disk.

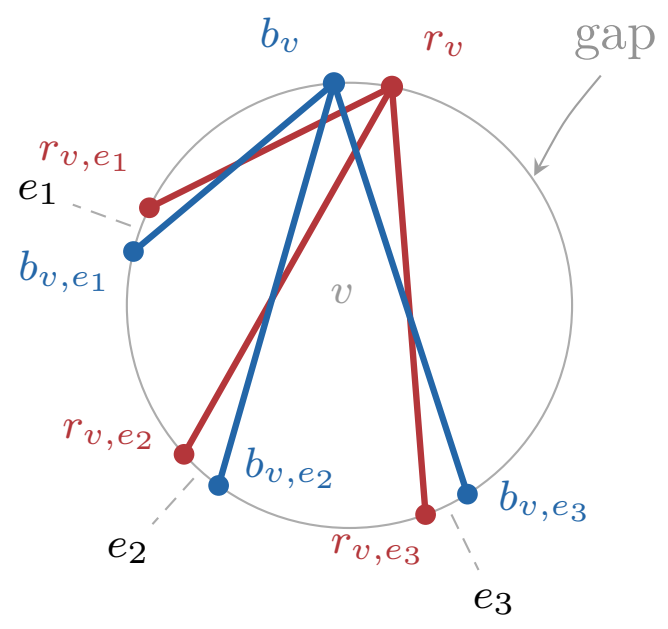
\begin{figure}[htb]
\centering
\begin{tikzpicture}[scale=.95,>=stealth]
  \definecolor{redpath}{RGB}{180,54,58}
  \definecolor{bluepath}{RGB}{35,102,168}

  \draw[gray!65] (0,0) circle (1.55);

  \coordinate (rv) at (79:1.55);
  \coordinate (bv) at (94:1.55);

  \foreach \i/\ra/\ba/\sp in {
    1/154/166/160,
    2/222/234/228,
    3/290/302/296
  }{
    \coordinate (r\i) at (\ra:1.55);
    \coordinate (b\i) at (\ba:1.55);

    \draw[gray!65,densely dashed]
      (\sp:1.62)--(\sp:2.05);
    \node[font=\scriptsize] at (\sp:2.30)
      {$e_{\i}$};

    \draw[redpath,very thick] (rv)--(r\i);
    \draw[bluepath,very thick] (bv)--(b\i);

    \fill[redpath] (r\i) circle (2.05pt);
    \fill[bluepath] (b\i) circle (2.05pt);
  }

  \fill[redpath] (rv) circle (2.2pt);
  \fill[bluepath] (bv) circle (2.2pt);

  \node[font=\scriptsize,above right=2pt,redpath]
    at (rv) {$r_v$};
  \node[font=\scriptsize,above left=2pt,bluepath]
    at (bv) {$b_v$};

  \node[font=\scriptsize,above left=2pt,redpath, xshift=3pt, yshift=1pt]
    at (r1) {$r_{v,e_1}$};
  \node[font=\scriptsize,below left=1pt,bluepath, xshift=+1pt, yshift=+5pt]
    at (b1) {$b_{v,e_1}$};


  \node[font=\scriptsize, redpath,  anchor=south east, xshift=-1pt, yshift=-3pt]
    at (r2) {$r_{v,e_2}$};
  \node[font=\scriptsize, bluepath, anchor=north west, xshift=1pt,  yshift=11pt]
    at (b2) {$b_{v,e_2}$};

  \node[font=\scriptsize, redpath,  anchor=north east, xshift=4pt, yshift=0pt]
    at (r3) {$r_{v,e_3}$};
  \node[font=\scriptsize, bluepath, anchor=south west, xshift=2pt,  yshift=-9pt]
    at (b3) {$b_{v,e_3}$};

  \node[font=\scriptsize,gray!80] at (-.05,.08)
    {$v$};

  \draw[gray!75,->]
    (1.75,1.60) .. controls (1.50,1.30) ..
    (35:1.57);
  \node[font=\small,gray!85] at (1.85,1.83)
    {gap};
\end{tikzpicture}
\caption{The red--blue resolution at a vertex.
The original spokes are shown dashed. Each color
connects its vertex in the chosen gap to the
corresponding ports on the disk boundary.}
\label{fig:vertexports}
\end{figure}

Outside the disks, run colored copies of each edge in a narrow corridor
along the original edge.  At a crossing of distinct original edges, the
red copy of one and the blue copy of the other cross once, with the
original sign.

Along a single edge, the red port is clockwise of the blue port at
\emph{both} ends.  Since the two outward directions are opposite, the
red and blue tracks lie on opposite relative sides of the edge at its
two ends.  We therefore join them with exactly one crossing in the
corridor, of sign $-1$ when the tracks have the same orientation
(Figure~\ref{fig:redblue}(b)).  No other corridor crossings are needed.

For edge chains $z,z'$, define
\[
 R'_D(z,z'):=I(z^{\mathrm{red}},(z')^{\mathrm{blue}}),
\]
where $z^{\mathrm{red}}$ follows $z$ along the red vertices and edges and $(z')^{\mathrm{blue}}$ follows $z'$ along the blue vertices and edges; this gives an 
unambiguous definition of a bilinear intersection number $R'_D$. We check that it satisfies~\eqref{eq:Rdirect}:
\begin{itemize}
\item For edges with no common endpoint, the red--blue crossings are
exactly the original crossings, with the same signs.
\item For distinct edges $[vp],[vq]$ directed away from $v$, the two
chords inside the disk cross if and only if $p\prec_v q$, and then
with sign $+1$ (Figure~\ref{fig:redblue}(a)).  All other crossings
contribute $I_D([vp],[vq])$.
\item For the same directed edge $[uv]$, the red and blue chords cross
with sign $+1$ inside each endpoint disk.  Their corridor crossing
has sign $-1$, giving $R'_D([uv],[uv])=1+1-1=1$.
\end{itemize}
Reversing an edge reverses the sign, so these cases determine the
pairing on all oriented edges and prove~\eqref{eq:Rdirect}.

Finally, the red and blue copies of any two chains are arbitrarily small
perturbations of the original chains.  If their endpoints are
non-interfering, such perturbations preserve their signed algebraic
intersection number.  Hence $R'_D(\gamma,\eta)=I_D(\gamma,\eta)$,
as claimed.
\end{proof}

Any two closed cycles in the plane have signed algebraic intersection number $0$ and
since they do not have interfering endpoints, the value of $R'_D$ must also be $0$.
\begin{corollary}\label{lem:R}
For any two cycles $z,z'$ in $G$,
\begin{equation}\label{eq:Rcycles}
R'_D(z,z')=0.
\end{equation}
\end{corollary}

An important property of the function $R'_D$ is the following. For a vertex $v$ and distinct neighbors $p,q$ note that 
$R'_D([vp],[vq]) + R'_D([vq],[vp]) = 1$. We let $R_D$ be $R'_D$ reduced modulo $2$. When working in a field extension of $\FF_2$ 
we will extend $R_D$ bilinearly to that field extension. Note that
\begin{equation}\label{eqtour}
R_D([vp],[vq]) + R_D([vq],[vp]) = 1,
\end{equation}
that is, $L(p,q):= R_D([vp],[vq])$ is a ``tournament'' (we will define this presently).

\section{Local detector}\label{sec:detector}

In this section we state the local detector theorem, a non-trivial algebraic fact which is at the core of the strengthened crossing-number argument. We postpone the proof to Sections~\ref{sec:universal}--\ref{sedif}.

We work over the rational-function field $\KK=\FF_2(u_1,\ldots,u_{N})$, for $N \geq 3$, where the $u_p$'s are independent variables and 
\begin{equation}\label{elam}
 \lambda_p:=\frac{1}{\prod_{q\ne p}(u_p+u_q)},
\end{equation}
are their Lagrange interpolation weights (these weights are rational functions of the $u_p$'s). Since the field has characteristic $2$
we have $X+Y=X-Y$.

A {\em tournament} is a matrix $L$ with entries in $\FF_2$ satisfying $L_{pq}+L_{qp}=1$ for $p\ne q$ (we don't put any restriction on the diagonal entries). Let 
\[
 s_1=\left\lceil\frac{N-3}{2}\right\rceil,
 \qquad s_2=\left\lfloor\frac{N-3}{2}\right\rfloor.
\]
Define the square $(s_1+1)(s_2+1)\times (s_1+1)(s_2+1)$ matrix $M_N(L)$, with rows and columns indexed by
$(a,b)$, where $0\le a\le s_1$ and $0\le b\le s_2$ as follows.
\begin{equation}\label{eq:M}
 M_N(L)_{(a,b),(c,d)}=
 \sum_{p\ne q}L_{pq}\lambda_p\lambda_q(u_p+u_q)
                       u_p^{a+d}u_q^{b+c}.
\end{equation}

\begin{theorem}\label{thm:detector}
For any $N\geq 3$ and for every tournament $L$ on $N$ vertices,
\[
\det M_N(L)\ne0\quad\text{in }
 \FF_2(u_1,\ldots,u_{N}).
\]
\end{theorem}

Let ${\cal P}_{a,b}$ be the space of polynomials in variables $X,Y$ with degree $\leq a$ in $X$ and degree $\leq b$ in $Y$.
We adopt the convention that \(\mathcal P_{a,b}=\{0\}\) when \(a<0\) or \(b<0\).
We will apply Theorem~\ref{thm:detector} in the following form.
\begin{corollary}\label{coap}
Suppose $P\in{\cal P}_{s_1,s_2}$ satisfies
\begin{equation}\label{eq:detector-functional}
 \sum_{p\ne q}L_{pq}\lambda_p\lambda_q(u_p+u_q)
 P(u_p,u_q) Q(u_p,u_q)=0
 \quad\text{for every }Q\in{\cal P}_{s_2,s_1}.
\end{equation}    
Then $P$ is the zero polynomial.
\end{corollary}

\section{The crossing-number lower bound}
\label{sec:odd}
Throughout this section $n=2s+1\ge5$ and $D$ is an arbitrary drawing of $K_n$.  Let $t_1,\dots,t_n$ be independent variables, one
for each vertex. We will work over a field $\KK\supseteq\FF_2(t_1,\ldots,t_n)$. Define the Lagrange coefficients
\[
 \omega_i=\frac1{\prod_{j\neq i}(t_i+t_j)}.
\]
Applying Lagrange interpolation to $x^a$ for $0\leq a\leq n-1$ at $t_1,\dots,t_n$ gives
\begin{equation}\label{elag}
x^a = \sum_{i=1}^n t_i^a \prod_{j\neq i}\frac{(t_j+x)}{(t_i+t_j)} = \sum_{i=1}^n \omega_i t_i^a \prod_{j\neq i} (t_j + x).
\end{equation}
Looking at the coefficient of $x^{n-1}$ on both sides we obtain
\begin{equation}\label{elag2}
\sum_i\omega_i t_i^a=0\ \mbox{for}\ a\in\{0,\dots,n-2\}\quad\mbox{and}\quad\sum_i\omega_i t_i^{n-1}=1.
\end{equation}

We say that a polynomial $P(x_1,\dots,x_n)$ is {\em symmetric} in $(x_i,x_j)$ if swapping 
$x_i$ and $x_j$ does not change its value, that is, 
\begin{equation}\label{syme}
P(x_1,\dots,x_i,\dots,x_j,\dots,x_n) = P(x_1,\dots,x_j,\dots,x_i,\dots,x_n).
\end{equation}
Let ${\cal P}_k$ be the family of $4$-variable polynomials that have degree at most $k$ in each variable.  Let $\SSS_k\subseteq {\cal P}_k$ contain polynomials 
that are symmetric in $(x,u)$ and in $(y,v)$, that is $P(u,x,y,v)=P(x,u,y,v)=P(x,u,v,y)$. The symmetry condition is equivalent to having the coefficient of 
$x^a u^b y^c v^d$ depend only on the two unordered exponent pairs $\{a,b\}$ and $\{c,d\}$. This implies 
\begin{equation}\label{eq:dimS}
\dim\SSS_k=\binom{k+2}{2}^{\!2}=Z(2k+5).
\end{equation}
The following definition is specific for fields of characteristic two. We say that a polynomial $P(x_1,\dots,x_n)$ is {\em alternating} in $(x_i,x_j)$  if it is symmetric in $(x_i,x_j)$ and for $x_i=x_j$ the value of the polynomial is zero, that is, 
\begin{equation}\label{alte}
P(x_1,\dots,x_i,\dots,x_i,\dots,x_n) = 0.
\end{equation}
Let ${\cal A}_k\subseteq\SSS_k$ contain polynomials that are alternating in $(x,u)$ and in $(y,v)$. The alternation
condition is equivalent to having the coefficient of $x^a u^b y^c v^d$ depend only on the two unordered exponent pairs $\{a,b\}$ and $\{c,d\}$, with the coefficient forced to zero 
if $a=b$ or $c=d$. This implies 
\begin{equation}\label{eq:dimSA}
\dim{\cal A}_k=\binom{k+1}{2}^{\!2}=Z(2k+3).
\end{equation}

From~\eqref{elag2} and Corollary~\ref{lem:R} we obtain the following.
\begin{lemma}\label{lem:tests}
For every polynomial $H(x,u,y,v)$ in ${\cal P}_{n-2}$
\begin{equation}\label{eq:fulltest}
 \sum_{i,j,k,\ell=1}^{n}\omega_i\omega_j\omega_k\omega_\ell
 R_D([ij],[k\ell])H(t_i,t_j,t_k,t_\ell)=0.
\end{equation}
For every polynomial $H(x,u,y,v)$ in ${\cal A}_{n-2}$
\begin{equation}\label{eq:wedge}
 \sum_{i<j}\sum_{k<\ell}\omega_i\omega_j\omega_k\omega_\ell
 R_D([ij],[k\ell])H(t_i,t_j,t_k,t_\ell)=0.
\end{equation}
\end{lemma}

\begin{proof}
Note that the sum in~\eqref{eq:fulltest} is effectively over $i\neq j$ and $k\neq\ell$ since, by convention, we let $[ii]=0$.
Let $0\le a,b\le n-2$. The boundary of the chain $z_{ab}=\sum_{i,j}\omega_i\omega_jt_i^at_j^b[ij]$ 
is 
\begin{equation}\label{chain}
\partial z_{ab}=\sum_{i,j}\omega_i\omega_jt_i^at_j^b([i]+[j]) = 
\sum_i [i]\omega_i t_i^{a} \sum_j \omega_j t_j^b 
+\sum_j [j]\omega_j t_j^{b} \sum_i \omega_i t_i^a = 0, 
\end{equation}
using~\eqref{elag2} in the final step.
Expanding $R_D(z_{ab},z_{cd})=0$ yields \eqref{eq:fulltest} for monomial $x^a u^b y^c z^d$. By linearity this extends to all $H$.

For the alternating case, using $[ij]=[ji]$ we can write $z_{ab}=\sum_{i<j}\omega_i\omega_j (t_i^at_j^b + t_i^b t_j^a) [ij]$.
Expanding \(R_D(z_{ab},z_{cd})=0\) yields \eqref{eq:wedge} for $(x^au^b+x^bu^a)(y^cv^d+y^dv^c)$. These polynomials span \(\mathcal A_{n-2}\), 
so the identity follows by linearity to all $H$.
\end{proof}

\begin{lemma}\label{lem:eval}
Suppose $P\in\SSS_{s-2}$ satisfies
\begin{equation}\label{hypoe}
 P(t_i,t_j,t_k,t_\ell)=0\quad\mbox{for all distinct}\ i,j,k,\ell\ \mbox{such that}\ R_{D}([ij],[k\ell])=1.
\end{equation}
Then $P$ is the zero polynomial.  
\end{lemma}

\begin{proof}
Assume $P(x,u,y,v)\in\SSS_{s-2}$ satisfies~\eqref{hypoe}, $P\neq 0$, and choose such $P$ with minimum degree in $x$. 
Note that for any four distinct
vertices $i,j,k,\ell$ we have
\begin{equation}\label{eeeer}
 R_D([ij],[k\ell])P(t_i,t_j,t_k,t_\ell)=0,
\end{equation}
because the assumption implies that at least one of the two terms in the product is $0$. 

Fix a vertex $i$. Let $u_1,\dots,u_{2s}$ be the elements of $\{t_1,\dots,t_n\}\setminus\{t_i\}$. Define
$$U_i(x)=\prod_{j\in\{1,\dots,s\}} \frac{x+u_j}{t_i+u_j}\quad\mbox{and}\quad V_i(x)=\prod_{j\in\{s+1,\dots,2s\}}\frac{x+u_j}{t_i+u_j}.$$ 
Note that $U_i,V_i$ are of degree $s$ and we have
\begin{equation}\label{eq:UiVi}
 U_i(t_j)V_i(t_j)=\begin{cases}1,&j=i,\\0,&j\ne i.\end{cases}
\end{equation}
Take any $Q\in\PP_{s-1,s-1}$ and let
\begin{equation}\label{eq:locpoly}
 H=P(x,u,y,v)(u+v)(u+y)(v+x)U_i(x)V_i(y)Q(u,v).
\end{equation}

Every variable in $H$ has degree at most $2s-1=n-2$ and hence we can apply~\eqref{eq:fulltest}. By equation~\eqref{eeeer} all independent-edge terms
are zero. The term $(u+v)(u+y)(v+x)$ zeroes out terms with $u=v$, $u=y$, and $v=x$. Hence the only terms that survive in the sum~\eqref{eq:fulltest} applied 
to $H$ are the ones with $x=y$ and with $x=y=t_i$ (the last equality is forced by~\eqref{eq:UiVi}).  Thus
\begin{equation*}
\begin{split}
 0 & =  \sum_{i',j,k,\ell\in [n]}\omega_{i'}\omega_j\omega_k\omega_\ell R_D([i'j],[k\ell])H(t_{i'},t_j,t_k,t_\ell) \\
  & = \omega_i^2 \sum_{\substack{j\neq\ell \in [n]\\ j,\ell\neq i}} \omega_j\omega_\ell R_D([ij],[i\ell]) (t_i+t_j) (t_i+t_\ell) (t_j+t_\ell) P(t_i,t_j,t_i,t_\ell) Q(t_j,t_\ell)\\
  & = \omega_i^2  \sum_{p\ne q\in [2s]} L(p,q)(u_p+u_q) \lambda_p\lambda_q P(t_i,u_p,t_i,u_q)Q(u_p,u_q), 
\end{split}
\end{equation*}
where $L$ is the tournament for vertex $i$ (see~\eqref{eqtour}) and in the last sum $p,q$ are such that $u_p=t_j, u_q=t_\ell$ and hence
\[
 \lambda_p= \frac1{\prod_{r\ne p}(u_p+u_r)}
 = \frac1{\prod_{k\neq i,j}(t_k+t_j)} = (t_i+t_j)\frac1{\prod_{k\neq j}(t_k+t_j)} =(t_i+t_j)\omega_j,
\]
and similarly $\lambda_q = (t_i+t_\ell)\omega_\ell$. Corollary~\ref{coap} yields
\begin{equation}\label{eq:localzero}
P(t_i,u,t_i,v)=0\quad\text{as a polynomial in }u,v.
\end{equation}

This holds for every vertex $i$.  Since $P(x,u,x,v)$ has degree at
most $2s-4<n$ in $x$, it follows that $P(x,u,x,v)\equiv0$.  Hence
$x+y$ divides $P(x,u,y,v)$.  By symmetry the product
\begin{equation}\label{eq:Cfactor}
 C=(x+y)(x+v)(u+y)(u+v)
\end{equation}
also divides $P$.  But $P/C\in\SSS_{s-2}$ satisfies the same vanishing
hypothesis, since $C$ is nonzero on four distinct vertex labels, and
$\deg_x(P/C)=\deg_xP-2$.  This contradicts the choice of $P$.
Therefore $P=0$.
\end{proof}

\begin{proof}[Proof of Theorem~\ref{thm:main}]
We can assume that $n \geq 5$. Let $w={\mathrm{iocr}(D)}$ be the number of independent pairs of edges crossing oddly in $D$ and define
\[
 \Theta=(x+u)(y+v)(x+y)(x+v)(u+y)(u+v).
\]
The purpose of the $(x+u)(y+v)$ terms is to make the polynomial alternating; 
the purpose of the remaining terms is to eliminate dependent edges. 

Let $P(x,u,y,v),Q(x,u,y,v)\in\SSS_{s-2}$. Then $PQ\Theta\in{\cal A}_{n-2}$. Let 
$T_D$ be the set of ordered pairs $((i,j),(k,\ell))$ where $i<j$, $k<\ell$, $i,j,k,\ell$ are all distinct  
and $R_D([ij],[k\ell])=1$. Note that $|T_D|=2w$, since each pair of independent edges crossing oddly occurs twice 
in $T_D$. From~\eqref{eq:wedge}, applied to $H=PQ\Theta\in{\cal A}_{n-2}$ we have
\begin{equation}\label{eq:wedgess}
 \sum_{((i,j),(k,\ell))\in T_D} 
 \alpha_{i,j,k,\ell}P(t_i,t_j,t_k,t_\ell)Q(t_i,t_j,t_k,t_\ell)=0, 
\end{equation}
 where
\begin{equation}\label{ealp}
\alpha_{i,j,k,\ell}=\omega_i\omega_j\omega_k\omega_\ell(t_i+t_j)(t_i+t_k)(t_i+t_\ell)(t_j+t_k)(t_j+t_\ell)(t_k+t_\ell).
\end{equation}
Note that for $((i,j),(k,\ell))\in T_D$ all factors on the right-hand-side of~\eqref{ealp}  are non-zero and hence $\alpha_{i,j,k,\ell}\neq 0$.
Also note that~\eqref{eq:wedgess} is true for every pair of 
polynomials $P(x,u,y,v),Q(x,u,y,v)\in\SSS_{s-2}$.

Let $d=\dim\SSS_{s-2}$ and suppose $w<d$. By Lemma~\ref{lem:eval}, evaluation on $T_D$ determines every polynomial in $\SSS_{s-2}$.  
Choose a subset $T_1\subseteq T_D$ of exactly $d$ points such that for every set of values on $T_1$ there is a unique polynomial
in $\SSS_{s-2}$ with those values (pick a basis $B_1,\dots,B_d$ of $\SSS_{s-2}$, for each $\alpha\in T_D$ consider the vector 
$(B_1(\alpha),\ldots,B_d(\alpha))$ and use the Steinitz lemma to pick $d$ linearly independent vectors). Let $T_2=T_D\setminus T_1$.  Since $|T_2|=2w-d<d$
there is a nonzero $Q(x,u,y,v)\in\SSS_{s-2}$ vanishing at every point of $T_2$. Hence for every $P(x,u,y,v)\in\SSS_{s-2}$ we have
\begin{equation}\label{eq:wedgess2}
 \sum_{((i,j),(k,\ell))\in T_1} 
 \alpha_{i,j,k,\ell}P(t_i,t_j,t_k,t_\ell)Q(t_i,t_j,t_k,t_\ell)=0.
\end{equation}
Since the values of $P$ are unrestricted (by our choice of $T_1$) we can conclude that $Q(t_i,t_j,t_k,t_\ell)=0$ for every $((i,j),(k,\ell))\in T_1$.
Thus $Q$ vanishes on all of $T_D$ (which is equivalent to \eqref{hypoe}) and contradicts Lemma~\ref{lem:eval}. Hence $w\ge d=\dim\SSS_{s-2}$.

This proves $\iocr(K_n) \geq  Z(n)$ for odd $n=2s+1$. The lower bound for even $n$ then follows a standard counting argument going back to Guy~\cite{G72}.
The upper bound is shown in the paper by Harary and Hill~\cite{HH63} that started it all.
\end{proof}

\section{The local detector: the algebraic argument} \label{sec:universal}

We now prove Theorem~\ref{thm:detector}. We fix a tournament $L$ and use definitions~\eqref{elam}
and~\eqref{eq:M}. Let $\VV=\PP_{s_1,s_2}$. The matrix $M_N$ represents the following bilinear form on $\VV$
\begin{equation}\label{eq:B}
 \BB(f,g)=\sum_{p\ne q}L_{pq}\lambda_p\lambda_q(u_p+u_q)
          f(u_p,u_q)g(u_q,u_p).
\end{equation}
We also use an auxiliary bilinear form, obtained by inserting one
more factor $u_p+u_q$:
\begin{equation}\label{eq:A}
 \Aux(f,g)=\sum_{p\ne q}L_{pq}\lambda_p\lambda_q(u_p+u_q)^2
          f(u_p,u_q)g(u_q,u_p).
\end{equation}
If $f,g,(X+Y)f\in\VV$ then
\begin{equation}\label{eq:relation}
\Aux(f,g)=\BB((X+Y)f,g).
\end{equation}
Let $\UU=\{g\in\PP_{s_1-1,s_1-1}:g(X,Y)=g(Y,X)\}\subseteq\VV$ and for $g\in\UU$ let 
\[
 \Phi(g)=\Aux(g,Xg).
\]
The proof needs the following four lemmas. Their proofs are given in
Sections~\ref{sec:derivative} and~\ref{sedif}.

\begin{lemma}\label{lem:diag}
The form $\BB$ is symmetric and, for every $f\in\VV$,
\begin{equation}\label{eq:diag}
\BB(f,f)=\sum_{p}
       u_p\bigl(\lambda_p f(u_p,u_p)\bigr)^2.
\end{equation}
We also have
\begin{equation}\label{eq:diagdivide}
 \BB(f,f)=0\quad\Longrightarrow\quad (X+Y)\mid f.
\end{equation}
\end{lemma}

\begin{lemma}\label{lem:half}
The form $\Aux$ is symmetric and for $g_0,g_1\in\UU$
\begin{equation}\label{eq:isotU}
 \Aux(g_0,g_1)=0.
\end{equation}
\end{lemma}

\begin{lemma}\label{lem:decom}
Every $h\in{\cal P}_{s_1-1,s_2-1}$ has a unique decomposition
\begin{equation}\label{eq:split}
h=h_0+Xh_1,\qquad h_0,h_1\in\UU.
\end{equation}
Explicitly,
\begin{equation}\label{eq:division}
 h_1=\frac{h(X,Y)+h(Y,X)}{X+Y},\qquad h_0=h+Xh_1.
\end{equation}
\end{lemma}

\begin{lemma}\label{lem:bridge}
For every non-zero $g\in\UU$ we have $\Phi(g)\neq 0$.
\end{lemma}

\begin{proof}[Proof of Theorem~\ref{thm:detector}.]
Suppose $f\in\VV$ satisfies $\BB(f,v)=0$ for every $v\in\VV$. We will argue that the only way this can happen is if $f=0$ (and hence $\BB$ is nonsingular). 
Taking $v=f$ and using Lemma~\ref{lem:diag} we obtain $f=(X+Y)h$, where $h\in{\cal P}_{s_1-1,s_2-1}$.  
By~\eqref{eq:relation} we have 
\begin{equation}\label{ahv}
\Aux(h,v)=0\ \mbox{for every}\ v\in\VV.
\end{equation}
Write $h=h_0+Xh_1$ with $h_0,h_1\in\UU$ using Lemma~\ref{lem:decom}.  Using 
Lemma~\ref{lem:half} we have
\[
 0=\Aux(h,h_1)=\Aux(h_0,h_1)+\Aux(Xh_1,h_1)
   =\Aux(h_1,Xh_1)=\Phi(h_1).
\]
Lemma~\ref{lem:bridge} now implies $h_1=0$ and hence $h=h_0\in\UU$. Taking $v=Xh_0$ in~\eqref{ahv} gives
$$0=\Aux(h_0,Xh_0)=\Phi(h_0),$$ 
and hence $h_0=0$. Thus $h=f=0$. We proved that $\BB$, and therefore $M_N(L)$, is nonsingular.
\end{proof}

\section{Basic properties of ${\cal A},{\cal B}$ - proofs of Lemmas~\ref{lem:diag},~\ref{lem:half},~\ref{lem:decom}}
\label{sec:derivative}

We use the following fact about  Lagrange interpolation coefficients (see argument after~\eqref{elag} for a proof). For $0\leq a\leq N-2$
\begin{equation}\label{eq:moments}
\sum_p\lambda_pu_p^a=0\quad\mbox{and}\quad \sum_p\lambda_pu_p^{N-1}=1.
\end{equation}
Note that~\eqref{eq:moments} implies that for a polynomial $f$ of degree $\leq N-2$ we have $\sum_p \lambda_p f(u_p)=0$. Similarly we can
conclude that for every $f(X,Y)$ that has degree $\leq N-2$ in $X$ or in $Y$
\begin{equation}\label{eq:moments2}
\sum_{p,q} \lambda_p\lambda_q  f(u_p,u_q) = 0.
\end{equation}
If $f\in{\cal P}_{N-2,N-2}$ is symmetric then the off-diagonal entries in~\eqref{eq:moments2}
can be paired up (since $\lambda_p\lambda_q  f(u_p,u_q) = \lambda_q\lambda_p  f(u_q,u_p)$) to cancel, hence
\begin{equation}\label{eq:moments3}
\sum_{p} \lambda_p^2  f(u_p,u_p) = 0.
\end{equation}

\begin{proof}[Proof of Lemma~\ref{lem:diag}]
Let $f, g\in\VV$. Then $f(X,Y)g(Y,X)\in\PP_{N-3,N-3}$. Swapping $p,q$ and using $L_{pq}+L_{qp}=1$, we obtain
\[
\begin{split}
 \BB(f,g)+\BB(g,f)&=
 \sum_{p\ne q}\lambda_p\lambda_q(u_p+u_q) f(u_p,u_q)g(u_q,u_p)\\
&=\sum_{p,q}\lambda_p\lambda_q(u_p+u_q) f(u_p,u_q)g(u_q,u_p)=0,
\end{split}
\]
using~\eqref{eq:moments3} to extend the summation from $p\neq q$ to all pairs $p,q$, and~\eqref{eq:moments2} in the last equality. This shows that $\BB$ is symmetric.  

Let $f\in\VV$. Swapping $p,q$ and using $L_{pq}+L_{qp}=1$, we obtain
\[
\begin{split}
 \BB(f,f)&=
 \sum_{p<q}\lambda_p\lambda_q(u_p+u_q) f(u_p,u_q)f(u_q,u_p)\\
 &=\sum_{p\neq  q}\lambda_p\lambda_q u_p f(u_p,u_q)f(u_q,u_p)\\
  &=\sum_{p}\lambda_p^2 u_p f(u_p,u_p)^2 + \sum_{p,q}\lambda_p\lambda_q u_p f(u_p,u_q)f(u_q,u_p)\\
  &=\sum_{p}\lambda_p^2 u_p f(u_p,u_p)^2,
\end{split}
\]
again using~\eqref{eq:moments2} for the last equality.

To prove the final claim of the lemma assume $\BB(f,f)=0$. Derivatives of squares vanish in characteristic two
and hence for every $j\in [N]$ we have
\begin{equation*}
\begin{split}
0 &= \frac{\partial}{\partial u_j} \BB(f,f) = \sum_{p} \frac{\partial}{\partial u_j} \lambda_p^2 u_p f(u_p,u_p)^2 = \\
& = \sum_{p} \big(   \frac{\partial}{\partial u_j}\big(\lambda_p f(u_p,u_p)\big)^2 u_p + \big(\lambda_p^2 f(u_p,u_p)^2\big) \frac{\partial}{\partial u_j} u_p \big)
= \lambda_j^2 f(u_j,u_j)^2.
\end{split}
\end{equation*}
Since the degree of $f(X,X)$ is at most $N-3<N$ and it is zero for $N$ values we must have $f(X,X)=0$ and hence $X+Y$ divides $f(X,Y)$
(since $f$ viewed as a polynomial in $X$ has root $Y$).
\end{proof}

\begin{proof}[Proof of Lemma~\ref{lem:half}]
Let $f, g\in\VV$. Swapping $p,q$ and using $L_{pq}+L_{qp}=1$, we obtain
\[
\begin{split}
 \Aux(f,g)+\Aux(g,f)&=
 \sum_{p\ne q}\lambda_p\lambda_q(u_p+u_q)^2 f(u_p,u_q)g(u_q,u_p)\\
&=\sum_{p,q}\lambda_p\lambda_q(u_p+u_q)^2 f(u_p,u_q)g(u_q,u_p)=0,
\end{split}
\]
where in the last equality we used~\eqref{eq:moments2} (expanding the square and dealing with $u_p^2$ and $u_q^2$ separately). Thus $\Aux$ is symmetric.  

Let $g_0,g_1\in\UU$. Then $g:=g_0g_1$ is also symmetric. Swapping $p,q$ and using $L_{pq}+L_{qp}=1$, we obtain
\[
\begin{split}
 \Aux(g_0,g_1)&=
 \sum_{p<q}\lambda_p\lambda_q(u_p+u_q)^2 g (u_p,u_q) \\
 &=\sum_{p\neq  q}\lambda_p\lambda_q u_p^2 g(u_p,u_q)\\
  &=\sum_{p}\lambda_p^2 u_p^2 g(u_p,u_p) + \sum_{p,q}\lambda_p\lambda_q u_p^2 g(u_p,u_q) = 0,
\end{split}
\]
where in the last equality we used~\eqref{eq:moments3}---applied to the symmetric polynomial
$XYg(X,Y)$---and~\eqref{eq:moments2}. 
\end{proof}

\begin{proof}[Proof of Lemma~\ref{lem:decom}]
For $s_2=0$ the statement is trivially true (since the space for $h$ is $\{0\}$); assume $s_2\geq 1$.
The numerator in~\eqref{eq:division} is $0$ when $X=Y$ and hence it is divisible by $X+Y$. The 
division decreases the degree in each variable (if non-zero). Thus 
$h_1$ is a polynomial in $\UU$. It is symmetric since its expression is
invariant under swapping $X$ and $Y$.

The polynomial $h_0=h+Xh_1$ is also symmetric:
\begin{equation*}
\begin{split}
h_0(X,Y)& =h(X,Y)+Xh_1(X,Y)=h(X,Y)+ Y h_1(X,Y) + (X+Y)h_1(X,Y)  \\
& = h(Y,X) + Y h_1(Y,X)  = h_0(Y,X).
\end{split}
\end{equation*}
The degree of $Y$ in $h_0$ is at most $\max\{s_2-1,s_1-2\}\leq s_1-1$ (since the bound is true for $h$ and $h_1$) and
by symmetry the same bound holds for the degree of $X$. Hence $h_0\in\UU$.

Suppose $h_0+Xh_1=0$ and $h_0$ and $h_1$ are symmetric. Swapping variables yields $h_0+Yh_1=0$ which combined
with $h_0+Xh_1=0$ gives $(X+Y)h_1=0$ and hence $h_1=0$ and $h_0=0$. This implies uniqueness for any $h\in {\cal P}_{s_1-1,s_2-1}$
(taking the difference of two purported solutions).
\end{proof}

\section{The differentiation lemma}\label{sedif}

For $0\leq a\leq N-1$ let
\begin{equation}\label{eq:D}
 \DD_a=\sum_{r=1}^{N}\lambda_r u_r^a\frac{\partial}{\partial u_r}.
\end{equation}

\begin{lemma}
For $0\leq a\leq N-1$ we have the following
\begin{equation}\label{eq:Dfacts}
 \DD_a u_p=\lambda_pu_p^a,
 \qquad
 \DD_a\lambda_p=
 \begin{cases}
 0,&a\text{ even},\\
 \lambda_p^2u_p^{a-1},&a\text{ odd}.
 \end{cases}
\end{equation}
\end{lemma}

\begin{proof}
The first equation in~\eqref{eq:Dfacts} is immediate from the definition. For the second we first collect basic facts on $(\partial/\partial u_r)$.
We have 
\begin{equation}\label{epar1}
\frac{\partial\lambda_p}{\partial u_r} = \frac{\lambda_p}{u_p+u_r}
\end{equation}
for $p\neq r$
and hence 
\begin{equation}\label{epar2}
\lambda_r\frac{\partial\lambda_p}{\partial u_r} = \lambda_p\frac{\partial\lambda_r}{\partial u_p}
\end{equation}
for any $p,r\in [N]$.
Using the product rule we also have
\begin{equation}\label{echain}
\frac{\partial (\lambda_r u_r^a)}{\partial u_p} = u_r^a \frac{\partial \lambda_r}{\partial u_p} + \mathbf1_{r=p}\, a \lambda_p u_p^{a-1}
\end{equation}
For the second equation in~\eqref{eq:Dfacts}, using~\eqref{epar2} and then~\eqref{echain}, we obtain
\begin{equation*}
\DD_a \lambda_p = \sum_{r=1}^{N} \lambda_r u_r^a \frac{\partial \lambda_p}{\partial u_r} = \lambda_p \sum_{r=1}^{N} u_r^a \frac{\partial \lambda_r}{\partial u_p} = \lambda_p \frac{\partial}{\partial u_p}\!\left(\sum_{r=1}^{N} \lambda_r u_r^a\right) + a \lambda_p^2 u_p^{a-1}.
\end{equation*}
Applying~\eqref{eq:moments} yields~\eqref{eq:Dfacts}.
\end{proof}

\begin{proof}[Proof of Lemma~\ref{lem:bridge}]
Suppose that $0\ne g\in\UU$ satisfies $\Phi(g)=0$, and choose
such a $g$ with minimum degree in $X$. Write
\[
 G_{pq}=(u_p+u_q)g(u_p,u_q).
\]
In terms of $G$ we have
\begin{equation}\label{eq:Q}
 \Phi(g)=\sum_{p\ne q}L_{pq}\lambda_p\lambda_q u_qG_{pq}^2.
\end{equation}

For $0\le a\le s_1$, apply $\DD_{2a}$ to $\Phi(g)$. Note that all $\lambda_p$ are annihilated using~\eqref{eq:Dfacts} and 
all squares are annihilated (since we work in characteristic $2$). The only term that we need to account for is $u_q$ 
and hence
\[
 \DD_{2a}\Phi(g)
 =\sum_{p\ne q}L_{pq}\lambda_p\lambda_q^2u_q^{2a}G_{pq}^2\\
 =\sum_p\lambda_p
    \left(\sum_{q\ne p}L_{pq}\lambda_qu_q^aG_{pq}\right)^2,
\]
where in the second equality we used $L_{pq}^2=L_{pq}$ (since $L_{pq}
\in\{0,1\}$).

We next apply $\DD_{2b+1}$, where $0\le b\le s_2$. Again the square is annihilated and the only term we need to account for is $\lambda_p$. 
By~\eqref{eq:Dfacts},
\begin{align}\label{eq:bridge}
 \DD_{2b+1}\DD_{2a}\Phi(g)
 &=\sum_p\lambda_p^2u_p^{2b}
    \left(\sum_{q\ne p}L_{pq}\lambda_qu_q^aG_{pq}\right)^2
    \notag\\
 &=\left(\sum_{p\ne q}L_{pq}\lambda_p\lambda_q
            u_p^bu_q^aG_{pq}\right)^2=\BB(g,X^aY^b)^2.
\end{align}
We assumed $\Phi(g)=0$ and hence $\BB(g,X^aY^b)=0$ for all $0\le a\le s_1$ and $0\le b\le s_2$. These monomials form a basis of $\VV$, so
\[
 \BB(g,v)=0\qquad\text{for every }v\in\VV.
\]
Since $g\in\UU\subseteq\VV$ we have $\BB(g,g)=0$ and, by Lemma~\ref{lem:diag}, $g=(X+Y)h$. Since $g$ is symmetric, so is $h$; moreover,
$h\in\UU$, $h\ne0$, and $\deg_X h<\deg_X g$.
Using~\eqref{eq:relation}, we obtain
\[
 \Aux(h,v)=\BB(g,v)=0\qquad\text{for every }v\in\VV.
\]
Taking $v=Xh\in\VV$ gives $\Phi(h)=\Aux(h,Xh)=0$,
contradicting the minimality of $\deg_X g$.
Thus $\Phi(g)\ne0$ for every nonzero $g\in\UU$.
\end{proof}

\section{The complete bipartite case}\label{sec:bipartite}

We apply the same techniques to complete bipartite graphs. We follow the structure of the proof in Section~\ref{sec:odd}, using similar ingredients: a cycle identity,
a localization lemma, and then a second application of the identity
to obtain the tight bound. 

We first prove the lower bound on $\iocr(K_{m,n})$ for odd $m=2r+1$ and $n=2s+1$, where $r,s\geq1$. Fix an arbitrary drawing $D$ of $K_{m,n}$.
Let the partitions of vertices be $\{1,\ldots,m\}$ and $\{1',\ldots, n'\}$.

Assign algebraically independent variables $t_1,\ldots,t_m$ to the vertices in $\{1,\dots,m\}$ and $u_1,\ldots,u_n$ to the vertices in 
$\{1',\dots,n'\}$. We work over
\(
 \KK=\FF_2(t_1,\ldots,t_m,u_1,\ldots,u_n)
\)
and use the Lagrange weights
\[
 \alpha_i=\frac1{\prod_{h\ne i}(t_i+t_h)},
 \qquad
 \beta_p=\frac1{\prod_{q\ne p}(u_p+u_q)}.
\]
As in~\eqref{elag2}, Lagrange interpolation gives
\begin{equation}\label{eq:bip-moments}
 \sum_i\alpha_i t_i^a=0\quad(0\leq a\leq m-2),
 \qquad
 \sum_p\beta_p u_p^b=0\quad(0\leq b\leq n-2).
\end{equation}
Let $\SSS_{r,s}$ denote the space of polynomials
$P(x,u,y,v)$ having degree at most $r$ in each of $x,y$
and at most $s$ in each of $u,v$. Note that
\begin{equation}\label{eq:bip-dimension}
 \dim\SSS_{r-1,s-1}=r^2s^2=Z(m,n).
\end{equation}

\begin{lemma}\label{lem:bip-tests}
For every polynomial $H(x,u,y,v)$ with degree at most
$m-2$ in each of $x,y$ and at most $n-2$ in each of $u,v$,
\begin{equation}\label{eq:bip-cycle}
 \sum_{i,k=1}^m\sum_{p,q=1}^n
   \alpha_i\alpha_k\beta_p\beta_q
   R_D([ip'],[kq'])H(t_i,u_p,t_k,u_q)=0.
\end{equation}
\end{lemma}

\begin{proof}
For $0\leq a\leq m-2$ and $0\leq b\leq n-2$, let
\[
 z_{ab}=\sum_{i=1}^m\sum_{p=1}^n
           \alpha_i\beta_{p} t_i^a u_{p}^b [ip'].
\]
Then $\partial z_{ab} = 0$: as in Equation~\eqref{chain} we argue that the coefficient of $[i]$ is
$\alpha_i t_i^a\sum_{p}\beta_{p}u_{p}^b=0$, and the coefficient of $[p']$ is
$\beta_pu_p^b\sum_i\alpha_i t_i^a=0$,
using~\eqref{eq:bip-moments}. 
By Corollary~\ref{lem:R}, $R_D(z_{ab},z_{cd})=0$.
Expanding this equality gives~\eqref{eq:bip-cycle} for
$H(x,u,y,v)=x^au^by^cv^d$, where
$0\leq a,c\leq m-2$ and $0\leq b,d\leq n-2$.
These monomials span $\SSS_{m-1,n-1}$, so the
identity follows for every $H$.
\end{proof}

\begin{lemma}\label{lem:bip-eval}
Suppose $P\in\SSS_{r-1,s-1}$ satisfies
\begin{equation}\label{eq:bip-vanish}
 P(t_i,u_p,t_k,u_q)=0
 \quad\text{whenever }i\ne k,\ p\ne q,\
              R_D([ip'],[kq'])=1.
\end{equation}
Then $P$ is the zero polynomial.
\end{lemma}

\begin{proof}
Suppose otherwise, and choose a nonzero $P\in\SSS_{r-1,s-1}$
satisfying~\eqref{eq:bip-vanish} with minimum degree in $x$.
For all $i\ne k$ and $p\ne q$ we then have
\begin{equation}\label{eq:bip-kernel}
 R_D([ip'],[kq'])P(t_i,u_p,t_k,u_q)=0.
\end{equation}
Fix $j\in[m]$. Partition $[m]\setminus\{j\}$ into two
sets $J_1,J_2$ of size $r$, and define
\[
 U_j(x)=\prod_{h\in J_1}\frac{x+t_h}{t_j+t_h},
 \qquad
 V_j(y)=\prod_{h\in J_2}\frac{y+t_h}{t_j+t_h}.
\]
Both polynomials have degree $r$, and
\begin{equation}\label{eq:bip-local-delta}
 U_j(t_i)V_j(t_i)=
 \begin{cases}
 1,&i=j,\\
 0,&i\ne j.
 \end{cases}
\end{equation}
For any $Q(u,v)\in\PP_{s-1,s-1}$, take
\begin{equation}\label{eq:bip-local-test}
 H(x,u,y,v)=P(x,u,y,v)(u+v)U_j(x)V_j(y)Q(u,v).
\end{equation}
The degree of $H$ in $x$ and $y$ is at most
$(r-1)+r=2r-1=m-2$, and in $u$ and $v$ it is at most
$(s-1)+1+(s-1)=2s-1=n-2$. We may therefore apply
Lemma~\ref{lem:bip-tests}.

Terms with $p=q$ vanish because of the factor $u+v$.
Terms with $i\ne k$ vanish by~\eqref{eq:bip-kernel}
(and the terms with $p=q$ have already been removed).
In the remaining terms $i=k$, so~\eqref{eq:bip-local-delta}
forces $i=k=j$. Hence,
\begin{equation}\label{eq:bip-local-sum}
 0=\alpha_j^2\sum_{p\ne q}
   L_{pq}\beta_p\beta_q(u_p+u_q)
   P(t_j,u_p,t_j,u_q)Q(u_p,u_q),
\end{equation}
where $L$ is the tournament for $j$, that is, $L_{pq}=R_D([jp'],[jq'])$, see~\eqref{eqtour}.

Corollary~\ref{coap}, with $N=n=2s+1$,
$s_1=s_2=s-1$, and $\lambda_p=\beta_p$,
applies to the polynomial
$F(u,v)=P(t_j,u,t_j,v)\in\PP_{s-1,s-1}$.
Since~\eqref{eq:bip-local-sum} holds for every
$Q\in\PP_{s-1,s-1}$, it follows that
\begin{equation}\label{eq:bip-diagonal-A}
 P(t_j,u,t_j,v)=0\qquad\text{as a polynomial in $u,v$.}
\end{equation}

This is true for every $j\in[m]$.
The polynomial $P(x,u,x,v)$ has degree at most $2r-2<m$ in $x$ and vanishes for $x=t_1,\ldots,t_m$. Hence it is identically zero, which implies that $x+y$ divides $P(x,u,y,v)$.
The quotient $P/(x+y)$ is nonzero, belongs to
$\SSS_{r-1,s-1}$, and still satisfies~\eqref{eq:bip-vanish},
since $t_i+t_k\ne0$ whenever $i\ne k$.
Its degree in $x$ is smaller, contradicting our
choice of $P$. Thus $P=0$.
\end{proof}

\begin{proof}[Proof of Theorem~\ref{thm:bipartite}]
Let $w=\iocr(D)$, and let $T_D$ be the set of
ordered quadruples $(i,p,k,q)$ for which $i\ne k$,
$p\ne q$, and $R_D([ip'],[kq'])=1$.
Each unordered pair of independent, odd-crossing edges
appears twice, and hence $|T_D|=2w$.

Take arbitrary $P,Q\in\SSS_{r-1,s-1}$ and let
\[
 H(x,u,y,v)=P(x,u,y,v)Q(x,u,y,v)(x+y)(u+v).
\]
The degree of $H$ in $x,y$ is at most
$2(r-1)+1=m-2$, and its degree in $u,v$
is at most $2(s-1)+1=n-2$.
The factors $x+y$ and $u+v$ eliminate all
pairs of edges with a common endpoint.
Applying Lemma~\ref{lem:bip-tests}, we obtain
\begin{equation}\label{eq:bip-pairing}
 \sum_{(i,p,k,q)\in T_D}
 \gamma_{ipkq}\,
 P(t_i,u_p,t_k,u_q)Q(t_i,u_p,t_k,u_q)=0,
\end{equation}
where
\[
 \gamma_{ipkq}=
 \alpha_i\alpha_k\beta_p\beta_q(t_i+t_k)(u_p+u_q).
\]
Every $\gamma_{ipkq}$ appearing in~\eqref{eq:bip-pairing}
is nonzero, since $i\ne k$ and $p\ne q$.
The identity holds for every pair of polynomials
$P,Q\in\SSS_{r-1,s-1}$.

Let $d=\dim\SSS_{r-1,s-1}=r^2s^2$.
Suppose $w<d$. By Lemma~\ref{lem:bip-eval}, evaluation
on $T_D$ determines every polynomial in $\SSS_{r-1,s-1}$.
Choose a subset $T_1\subseteq T_D$ of exactly $d$ points
whose evaluation rows are linearly independent. Then
for any prescribed values on $T_1$ there is a unique
polynomial $P\in\SSS_{r-1,s-1}$ taking those values.

Let $T_2=T_D\setminus T_1$. As
$|T_2|=2w-d<d$, there exists a nonzero
$Q\in\SSS_{r-1,s-1}$ vanishing at every point of $T_2$.
For this choice of $Q$, equation~\eqref{eq:bip-pairing}
reduces, for every $P\in\SSS_{r-1,s-1}$, to
\begin{equation}\label{eq:bip-final}
 \sum_{(i,p,k,q)\in T_1}
 \gamma_{ipkq}\,
 P(t_i,u_p,t_k,u_q)Q(t_i,u_p,t_k,u_q)=0.
\end{equation}
Since the values of $P$ are unrestricted we can conclude 
$Q(t_i,u_p,t_k,u_q)=0$ for every point of $T_1$. Thus $Q$ vanishes
on all of $T_D$, contradicting Lemma~\ref{lem:bip-eval}.
Hence 
\(
 w\geq d=r^2s^2=Z(2r+1,2s+1).
\)

This proves $\iocr(K_{m,n}) \geq Z(m,n)$ in case both $m$ and $n$ are odd.
The remaining cases, where at least one of $m$ or $n$ is even, then follow from a standard counting argument which can be found in Guy~\cite{G68}.
%
%
Finally, the classical Zarankiewicz drawing
has exactly $Z(m,n)$ crossings
~\cite{Z53,openaiBipartite2026}.
Since $\iocr\leq\cro$, this gives the matching
upper bound and completes the proof.
\end{proof}

\section{Open Questions}

\paragraph{Harborth's Conjecture}

The best known upper bound $Z(n_1,\ldots,n_k)$ on the crossing number of the complete $k$-partite graph $K_{n_1,\ldots,n_k}$ was found by
Harborth~\cite{H71}, also see~\cite[Theorem 2.1]{CHN20}. Harborth conjectured that this bound is optimal. Can the edge chain approach be adapted to prove Harborth's conjecture? The case $k=3$ would be a good starting point.

\paragraph{How close are $\iocr$ and $\cro$?}

All the crossing number variants we have mentioned in this paper are polynomially related; this follows from $\cro(G) \leq \binom{2\iocr(G)}{2}$~\cite{PSS10}. While we do know that some of these crossing numbers differ, e.g.\ $\iocr$ and $\ocr$ as mentioned earlier, this question remains open for many other variants, e.g.\ for $\cro_-$ vs $\cro$ or, famously, $\pcr$ vs $\cro$. Maybe Tutte's crossing chains as developed in the OpenAI approach can lead to new insights on these questions. 

\bibliographystyle{plainurl}
\bibliography{iocr}

\end{document}